\documentclass{tMAA2e}

\usepackage{pgf,tikz}
\usetikzlibrary{calc}
\newcommand{\radius}{1.5cm}
\newcommand{\h}{4.71cm} 

\begin{document}


\title{Sphericons and D-forms: a crocheted connection}

\author{Katherine A.~Seaton$^{a}$$^{\ast}$\thanks{$^\ast$ Email: k.seaton@latrobe.edu.au
\vspace{6pt}} \\\vspace{6pt} $^{a}${\em{Department of Mathematics and Statistics, La Trobe University VIC 3086, Australia}}
 }

\maketitle

\begin{abstract}
Sphericons and D-forms are 3D objects created and described by artists, which have separately received attention in the mathematical literature in the last 15 or so years. The attempt to classify a seamed, crocheted form geometrically led to the observation, which appears not to have been previously made explicit, that these objects are related. General results concerning $(N,k)$-icons and seam-, D- and pita-forms are given. Instructions to  crochet  such forms are provided
 in the Appendix.


\begin{classcode} 51M04; 97M80 \end{classcode}

\end{abstract}

\section{Introduction}

In January 2016, a crocheter on the social networking site Ravelry \cite{rav} posted images of a coin purse she had made as a Christmas present, and started a discussion thread asking readers of a particular forum if they could help her identify the name of its 3D shape (hereafter termed {\em The Shape}).

\begin{figure}
\begin{center}
{ \resizebox*{8cm}{!}{\includegraphics{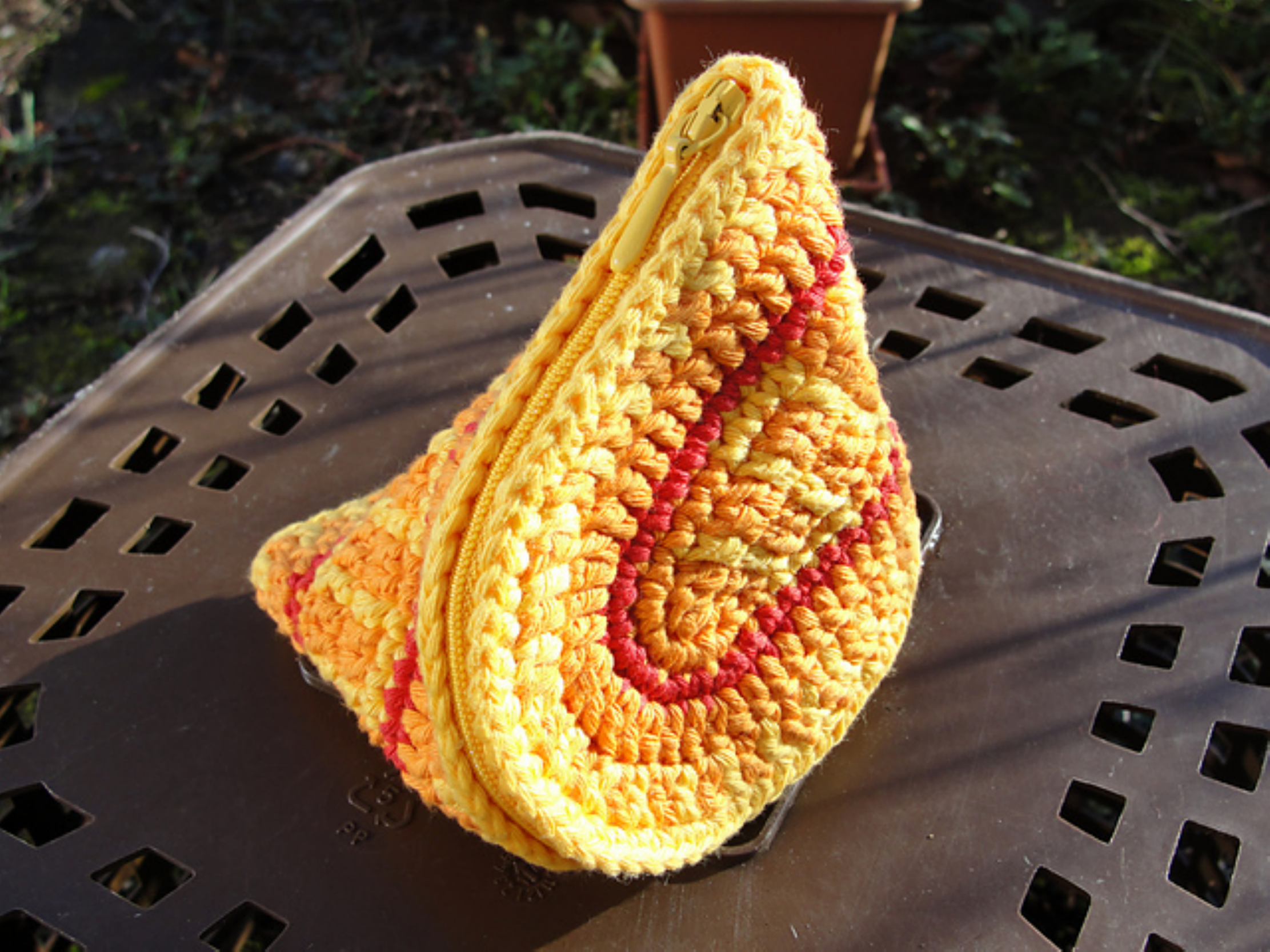}}}
\caption{\label{purse} Photograph of the coin purse, The Shape, made from a crocheted stadium. Courtesy of Annagret.}
\end{center}
\end{figure}

The purse, shown in Figure \ref{purse}, had been constructed following a pattern \cite{CB} by making a flat oval, and attaching a zipper to its edge in such a way that when the zipper closed, the purse was not folded over exactly in half onto itself (flat), but sat up, enclosing an intriguingly-shaped volume. The crocheter had read the German version \cite{PS} of an Ian Stewart column on the sphericon \cite{St}, and believed that The Shape might be sort-of-but-not-quite a sphericon. However, the paper templates for making one's own sphericons that she had seen \cite {PS, St} have curved edges only, which her oval did not.

\begin{figure}
\begin{center}
\begin{tikzpicture}
  \coordinate (A) at (0,0);
  \coordinate (B) at (0,2*\radius);
  \coordinate (C) at (\h,0);
  \coordinate (D) at (\h,2*\radius);
  \draw (A) -- node[above] {$h$} (C);
  \draw (C) arc [start angle=-90, end angle=90, radius=\radius] --
    (B) arc [start angle=90, end angle=270, radius=\radius];
  \draw[densely dashed] (A) -- node[left] {$2r$} (B);
  \draw[densely dashed] (C) -- node[right, pos=0.75] {$r$} (D);
  \node[circle, fill=black, minimum size=3pt, inner sep=0pt] at ($0.5*(C)+0.5*(D)$) {};
\end{tikzpicture}
\caption{\label{stad} A stadium, a rectangle of length $h$ inserted between two halves of a circle of radius $r$.}
\end{center}
\end{figure}
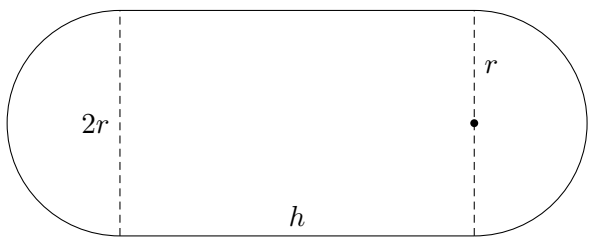

In fact  {\em oval} is rather  loose a term, and broadly means {\em egg-like}. The purse in question had been constructed from a {\em stadium}, the 2D shape made from a rectangle and two semi-circles, as shown in Figure \ref{stad}.

Over the next few weeks, various Ravelry members, including the author of this paper, joined in on-line discussion and brain-storming about The Shape. More than 500 of Ravelry's (admittedly 6 million) members read at least part of the thread (before, as threads do, it wandered off to shapes more generally). Patterns for objects with similar construction were pointed out \cite{ZS}. As well as the sphericon family, D-forms \cite{DP, Sh2} were mentioned.

So, was The Shape one of these, and if so, which was it? This seemingly simple question, prompted by a hand-made Christmas gift, has triggered this survey and
investigation of sphericons and D-forms. The main conclusions of this investigation are to be found in Section \ref{concl}.

The sphericon has been described as `[a] solid ... not as widely known as it should be' \cite{EW}, and conversations with mathematical colleagues, and presentations about this work,
have shown this to be so. Thus, an introduction to the history,  properties and generalisations of sphericons, and of D-forms, make up the next two sections of the paper.
The classification of The Shape is then clear. The paper concludes with some more general results, summarised in Table ~\ref{classify}, and includes  templates
for the surfaces required to make several particular generalised
sphericons by D-forming,  thereby stitching (literally) a connection between the sphericon family and D-forms.

 \section{Sphericons: making things out of solids} \label{solid}
The sphericon, at least as a mathematical object, came to public attention in Stewart's Mathematical Recreations column \cite{St} in 1999. Colin Roberts, a joiner who had enjoyed geometry at school, sent Stewart his creation, after seeing him on television. The column was written shortly after this. Roberts had first made a sphericon from mahogany some thirty years earlier in 1967 \cite{PR}. Independently, in the early 1980s, an Israeli designer David Hirsch submitted a patent application for `meander motion' toys constructed on the same principles \cite{DH}. Roberts' motivation, however, lay in producing a one-surfaced object, with no hole (unlike the well-known M\"obius strip).

\begin{figure}
\begin{center}
{ \resizebox*{12cm}{!}{\includegraphics{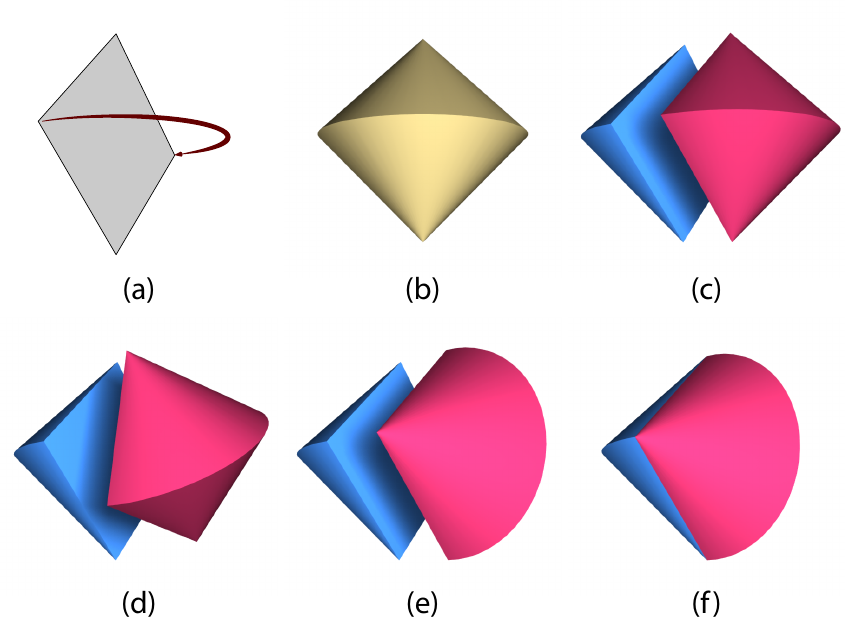}}}
\caption{\label{cut} Diagram showing the production of a sphericon from the solid of revolution of a square. Courtesy of Craig S. Kaplan.}
\end{center}
\end{figure}

\begin{figure}
\begin{center}
{ \resizebox*{6cm}{!}{\includegraphics{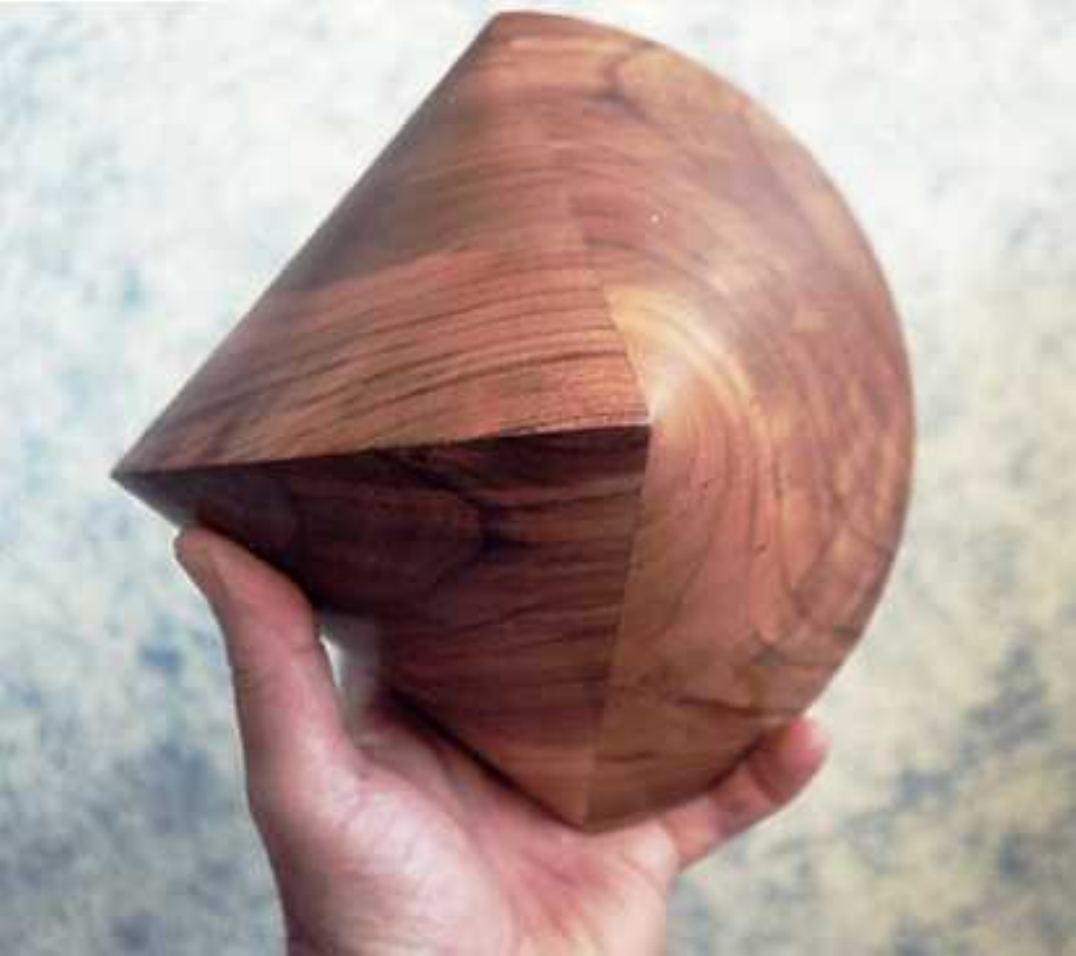}}}
\caption{\label{sphericon} A wooden sphericon. Permission given to reproduce from \cite{SS}.}
\end{center}
\end{figure}

The sphericon can be constructed from a pair of identical cones of apex angle 90 degrees joined at their bases (a {\em bicone}). Cut the bicone in half  through the two apices so as to reveal its square cross-section, twist one half through 90 degrees, and then `glue' the square faces together again, as shown in Figure \ref{cut}.  Especially when made of grained wood, it is a beautiful object.  It has a single, sinuous surface \cite{SS}, and two disjoint curved edges, set at right angles; see Figure~\ref{sphericon}. The surface has the property that if one traces a finger along the middle of it, one comes back to the starting point \cite{PR}. (Rather inaccurately from the point of view of mathematical terminology, this has been called a `continuous' surface \cite{PR, RK}; we will avoid that term.)
The two curved edges each end in vertices (four in total).

The sphericon rolls in an amusing manner, what Stewart terms `a controlled wiggle' \cite{St}, the meander motion of Hirsch's toys \cite{DH}. A {\em sphere} will roll in a straight line, and a {\em cone} rolls in a circle. The sphericon changes direction as the parts of its surface come into contact with the surface it is rolling on, but on average the motion is in a line. Hence the name {\em sphericon}. Sphericons are one of the intriguing solids that can be assembled, and their rolling behaviour predicted, at MoMath in New York \cite{MF}.

The sphericon featured next in the monthly Feature Column of the American Mathematical Society \cite{TP}. The sphericon was used as an interesting examplar of how  edges and  apices have to be treated carefully, in order to calculate the total Gaussian curvature of the figure to be the expected $4 \upi$.  Phillips introduced the \textit{hexasphericon} and \textit{octasphericon} and found connections between such higher order $N$-icons (to be defined below) and his work on mazes \cite{AT2}. Some particular closed space curves were found to lie on the sphericon surface  \cite{KM}, and an explicit hybrid function form (in cartesian coordinates) is given.

Perhaps not surprisingly, given its origins, wood-turners began elaborating upon the sphericon. A special join, called a paper join, enables the cutting-open process to be done without damage. The bicone can be replaced by any solid of revolution with appropriate cross-section symmetry to create new objects. Some makers introduced different profiles (producing what they termed {\em femispheres}) or non-convex cross-sections (the {\em streptohedrons}) \cite {Sp, Sp2}. Other makers \cite{SS} generated some of the convex geometric objects that are classified as the series of $N$-icons \cite{RK}. Sphericons (and these other elaborations) have been made of ceramic, glass, plastic (3D printing) and as metal sculptures, as well as of wood.

Think of the bicone as the solid of revolution of a square about its diagonal (as in the image (a) of Figure \ref{cut}). Now replace the square by any other $N$-vertex regular convex polygon, and use one of the
symmetry axes to generate a solid of revolution \cite{RK}. When $N$ is odd, this axis must run through one vertex and the midpoint of the opposing edge of the polygon. When $N$ is even, the
polygonal cross-section can be rotated about an axis through two opposite vertices or through the midpoints of two opposing edges.  For $N>4$, different `$N$-icons' can be produced by more than one rotation (in the cut-rotate-reglue process), due to the rotational symmetries of the cross-section; these may be completely different or may come in right- and left-hand pairs. A full classification of the family of solids requires one to specify $k \in \mathbb N$ where $\frac{2k \upi}{N}$ is the rotation (in radians) between cutting and regluing, termed the twist \cite{RK, PR}. It is also necessary to specify if the symmetry axis used in construction is vertex-to-vertex (vv) or midpoint-to-midpoint (mm) when $N$ is even. The notation $(N,k)$-icon will be used and, when $N$ is even, vv  or mm will be specified.

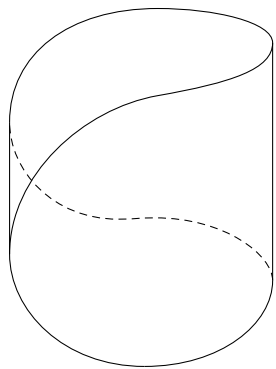
\begin{figure}
\begin{center}
\begin{tikzpicture}
  \coordinate (A) at (60pt,7.5pt);
  \coordinate (B) at (108.5pt,40pt);
  \coordinate (C1) at (64pt,63.75pt);
  \coordinate (C2) at (50pt,63.25pt);
  \coordinate (D) at (9pt,100pt);
  \coordinate (E) at (65pt,143pt);
  \coordinate (F) at (108.5pt,130pt);
  \coordinate (G) at (65pt,110pt);
  \coordinate (H) at (9pt,50pt);
  \draw (F) -- (B);
  \draw (D) -- (H);
  \draw (A) .. controls ($0.4*(A)+0.6*(B |- A)$) and ($0.5*(B)+0.5*(B |- A)$) .. (B);
  \draw[densely dashed]
    (B) .. controls ($0.6*(B)+0.4*(B |- C1)$) and ($0.4*(C1)+0.6*(B |- C1)$) ..
    (C1) .. controls ($0.5*(C1)+0.5*(C2 |- C1)$) and ($0.5*(C2)+0.5*(C1 |- C2)$) ..
    (C2) .. controls ($0.4*(C2)+0.6*(D |- C2)$) and ($0.5*(D)+0.5*(D |- C2)$) .. (D);
  \draw (D) .. controls ($0.4*(D)+0.6*(D |- E)$) and ($0.4*(E)+0.6*(D |- E)$) ..
    (E) .. controls ($0.5*(E)+0.5*(F |- E)$) and ($0.4*(F)+0.6*(F |- E)$) ..
    (F) .. controls ($0.5*(F)+0.5*(F |- G)$) and ($(G)+(10:1)$) ..
    (G) .. controls ($(G)+(190:1)$) and ($0.5*(H)+0.5*(H |- G)$) ..
    (H) .. controls ($0.5*(H)+0.5*(H |- A)$) and ($0.4*(A)+0.6*(H |- A)$) .. (A);
\end{tikzpicture}
\caption{\label{dual} The `dual sphericon', more precisely termed the mm $(4,1)$-icon.}
\end{center}
\end{figure}

Two special cases from the earlier literature will now be considered:

\textbf{The `dual' sphericon} Rotating a square about an axis through the midpoints of opposite edges results in a cylinder of equal height and diameter. Applying the cutting and rotating
process, an object with {\em two surfaces} and {\em one closed edge} is produced, as illustrated in Figure \ref{dual}. In the terminology just introduced, we would call it the mm $(4,1)$-icon.
This object is dual to the original sphericon (the vv $(4,1)$-icon) in a particular sense described by Kaufman \cite{RK}; the sphericon has {\em two disjoint edges} ending in vertices and {\em one surface} that can be traced as described
above. Springett \cite{Sp} described this object, which wood-turners term {\em side-cut},  as being `rather unexciting'\footnote{It will become apparent that this author does not share this opinion.}. In the cutting process, four semicircles are produced, which are then joined smoothly to the sides of the cylinder. The paper template for this object consists of not one stadium, but two. The mm $(4,1)$-icon is a possible two-sided die. Indeed, such dice --- and various $(N, k)$-icons --- are available from various 3D printing artists.

\begin{figure}
\begin{center}
{ \resizebox*{6cm}{!}{\includegraphics{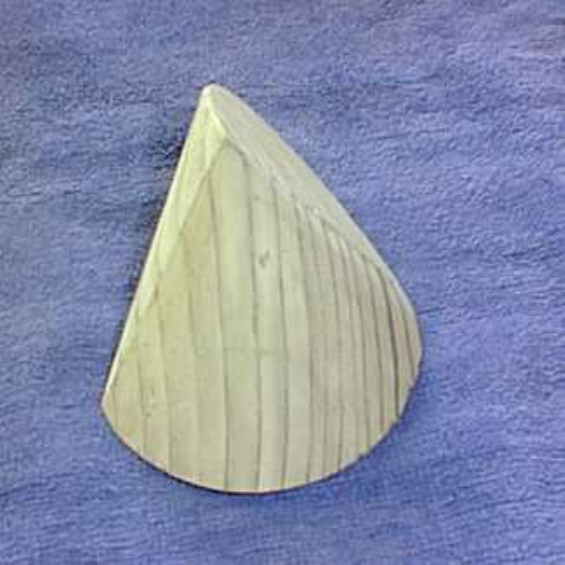}}}
\caption{\label{trisph} A wooden trisphericon, the $(3,1)$-icon. Permission given to reproduce from \cite{SS}. }
\end{center}
\end{figure}

\begin{figure}
\begin{center}
\begin{tikzpicture}
  \coordinate (A) at (0,0);
  \coordinate (B) at (0,2*\radius);
  \coordinate (C) at (2*\radius,2*\radius);
  \coordinate (D) at (2*\radius,4*\radius);
  \draw (A) -- node[left] {$2r$} (B);
  \draw (B) -- node[above] {$2r$} (C);
  \draw (C) -- node[right] {$2r$} (D);
  \draw[densely dashed] (A) -- (C);
  \draw[densely dashed] (B) -- (D);
  \draw (D) arc [start angle=90, end angle=180, radius=2*\radius]
    arc [start angle=90, end angle=270, radius=\radius]
    arc [start angle=270, end angle=360, radius=2*\radius]
    arc [start angle=-90, end angle=90, radius=\radius];
  \draw ($(B)+(0,-0.2*\radius)$) -- ++(0.2*\radius,0) -- ++(0,0.2*\radius);
  \draw ($(C)+(-0.2*\radius,0)$) -- ++(0,0.2*\radius) -- ++(0.2*\radius,0);
\end{tikzpicture}
\caption{\label{temp} Template for a trisphericon, the $(3,1)$-icon. It is formed from two quarter-circles of radius $2r$ and two semi-circles, of radius $r$.}
\end{center}
\end{figure}
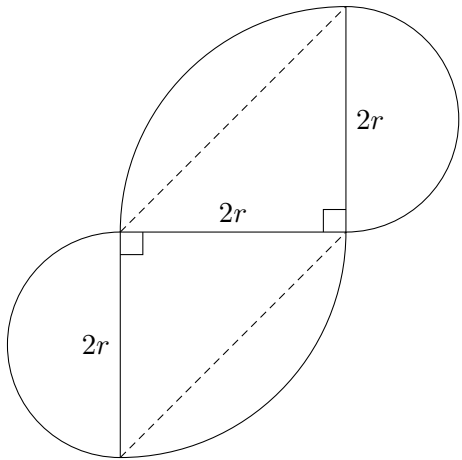

\textbf{The trisphericon} A {\em single} cone may be created as the solid of revolution of an equilateral triangle, the $N=3$ case. It can be cut to produce a triangular cross-section.
When one half is rotated through 120 degrees and the form reglued, the $(3,1)$-icon is formed. This is a solid with one curved surface, one \textsf{S}-shaped edge and 2 vertices; see Figure~\ref{trisph}. In the (Kaufman) sense described above, this object is self-dual. Following a similar argument to Stewart \cite{St} we can establish what the paper template for this object is (Figure \ref{temp}). This object has also been called  the `kleines sphericon' \cite{PS}, the `conicon' \cite{RK} or the `equilaticon' \cite{Charlie}. While `one-sided dice' based on the M\"obius band are available to purchase,  $(3,1)$-icons are even better candidates for this title, as they have the desirable feature of having no hole.

More generally, the vv $N$-icons, and the odd $N$ solids, can always be considered as formed by cutting, twisting and reglueing a solid made of $\lfloor \frac{N}{2}\rfloor$ sections of cones
(possibly including a cone with apex at infinity, a cylinder). Stadia may be present in the mm $N$-icons, as we can see from the four semicircles created in the cutting process.
New features arise as $N$ increases, such as surfaces enclosing other surfaces \cite{RK} ---  think of a velodrome, with the sloped surface around the centre.

Kaufman \cite{RK} also discusses a hybrid object, formed at the joining step from half of a vv $N$-icon and half of a mm $N$-icon (with the same value of $N$); we do not consider such objects in this paper.

It is worthwhile, at this point, if you have resisted the temptation so far,  to make a sphericon or trisphericon from paper or light card using templates such as those of P\"{o}ppe \cite{PS}. (Indeed, the crocheter mentioned in the Introduction had done exactly this with her family, which is why she was familiar with them.) Of course, you will not have a solid; you will have the surface and will manipulate and tape it to enclose a volume shaped like the solid.

\section{D-forms: making things out of sheets}
The D-forms entered mathematics in a similar way to the sphericon family--- but this time from design. They are 3D objects formed when the boundaries of two 2D shapes of equal perimeter length are joined.

\begin{figure}
\begin{center}
{ \resizebox*{6cm}{!}{\includegraphics{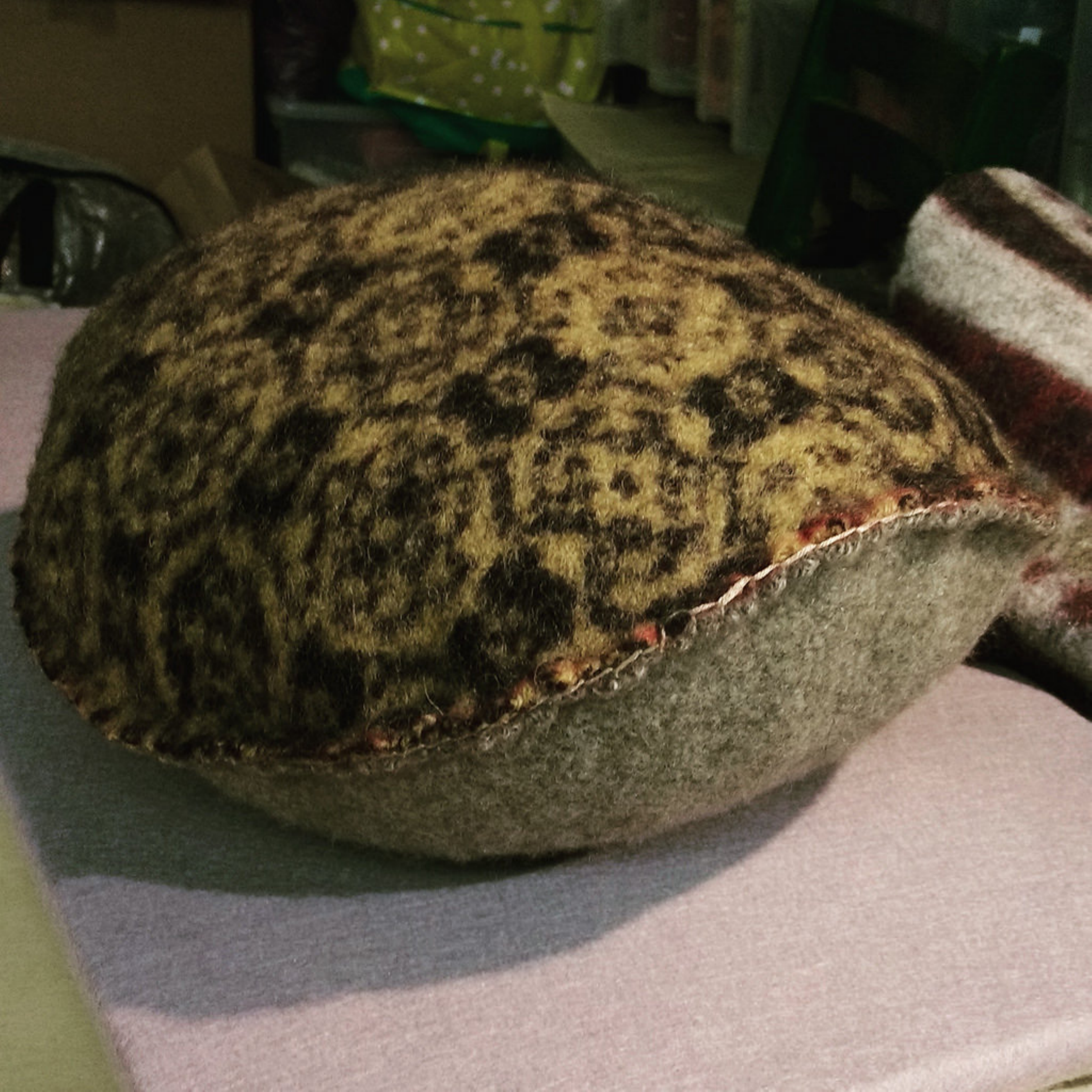}}}
\caption{\label{MS} A D-form made from two ellipses. Courtesy of Madeleine Shepherd.}
\end{center}
\end{figure}

The British designer and inventor Tony Wills says they came to him in a dream, as a way of expanding his working vocabulary of forms \cite{Wi}.
Wills has used D-forms in street architecture (benches, bins, bollards) \cite{WW}.  The constituent 2D shapes need to be made of a material that does not stretch or shear;
a typical cushion is not a D-form! Two circles cannot produce a (non-trivial) D-form, but two identical offset ellipses do, as seen in Figure ~\ref{MS}. The  cross-section of the object produced using two ellipses with their major axes set at right angles looks like the letter capital \textsf{D}. (D is for dream, for the shape of this cross-section and also for developable surfaces \cite{Sh2}.)

 As the shapes are joined, they inform each other as to how their joined boundary should curve in space. Wills observed that choosing different  points on the boundaries of the 2D shapes to start the seam between them changes the resulting D-form \cite{Wi, PB}. He  investigated joining the boundary of a circle to that of a square (i.e. not a smooth boundary) producing an object he termed a Squaricle, and also used non-convex shapes. Baseball construction \cite{BB} is an example of D-forming that long pre-dates Wills; two flat (non-convex) pieces of leather have been designed (by trial and error initially) to be joined in such a way that their seam lies on a sphere. The initial 3D form is not spherical --- the leather is subsequently stretched over a hard spherical core.

D-forms are piece-wise {\em developable surfaces}, in that they can be cut open and each piece flattened, without creasing, to be planar \cite{Sh1, PW}, although they have been described thus far as being  formed by the reverse process.  Developable surfaces have the advantage for designers, artists and builders that they can be constructed from rolled sheet metal; even ships and airplanes are made this way \cite{AKA}.  D-forms have been observed to have much in common with various sculptural works \cite{Sh1}. Predicting in advance what the D-form will look like, and hence what 2D shapes to use, is not easy, though computational methods (based on meshes) have been developed \cite{GAS}. A recent paper \cite{Bridgesnew} describes software used (in the education of architecture students) to match the perimeter length of different planar shapes, scaling them to be isoperimetric, in order that they can be used to create D-forms.

As with the sphericon, making a paper model, right now, will help you understand how D-forms evolve in front of you and how beautiful they can be. You can use a template \cite{LB, PB}  --- or you can easily cut your own using a double layer of paper.

 D-forms were introduced to the mathematical world by Wills and Sharp before their publications \cite{Sh1,Sh2}. Questions  \cite[][p. 418]{PW} regarding their creasing, and whether or not a D-form is always the convex hull of a space curve, were posed in 2001 and addressed in 2007. It is important to note that these discussions restricted the boundaries of the 2D shapes to being convex \cite{DO}; thus the D-forms of the mathematical literature are a particular subset of the D-forms of the design literature.

Demaine and O'Rourke also introduced the concept of a {\em pita-form} \cite{DO}. A pita-form is defined to be the 3D shape obtained by glueing (or zipping) the boundary of a {\em single} convex 2D shape to itself. Demaine returned to the open questions surrounding  D-forms with a different co-author in 2010 \cite{DP}. They established that both D-forms and pita-forms are the convex hull of their seams, that D-forms are always crease-free away from the seam, but that a pita-form may have at most one crease, between the two endpoints of its seam. (Of course, `crease' has to be defined a bit more carefully (in terms of differentiability properties) than by simply looking at a paper model, as this initially led to the proposition that pita-forms also might never have creases \cite{DO}.) Also introduced were {\em seam forms}, a natural generalization to more flat (convex) pieces than two.

(It may be pertinent to mention one more concept introduced by Demaine and collaborators, that of {\em zipper unfolding} of a polyhedron \cite{DD}, to prevent confusion. In this case, the surface being zipped is deliberately folded, marking it into the polygons that constitute the faces of the polyhedron.)

\section{Stitching it up: making things out of fibre}

As mentioned in the Introduction, both sphericons and D-forms were suggested as describing The Shape (Figure~\ref{purse}), formed from a crocheted stadium, lined so that there is some rigidity to the surface and no shearing or stretching, before the zipper is added.

The astute reader will by now have realised that, simply, it is a pita-form, made from a stadium.

This author has experimented with making various pita-forms from crocheted stadia, varying the ratio between the rectangle's height and the radius of the semi-circles, and with different starting points for the seam (Figure~\ref{pita}). Tips for making your own from fibre are given in Appendix A.

\begin{figure}
\begin{center}
{ \resizebox*{17cm}{!}{\includegraphics{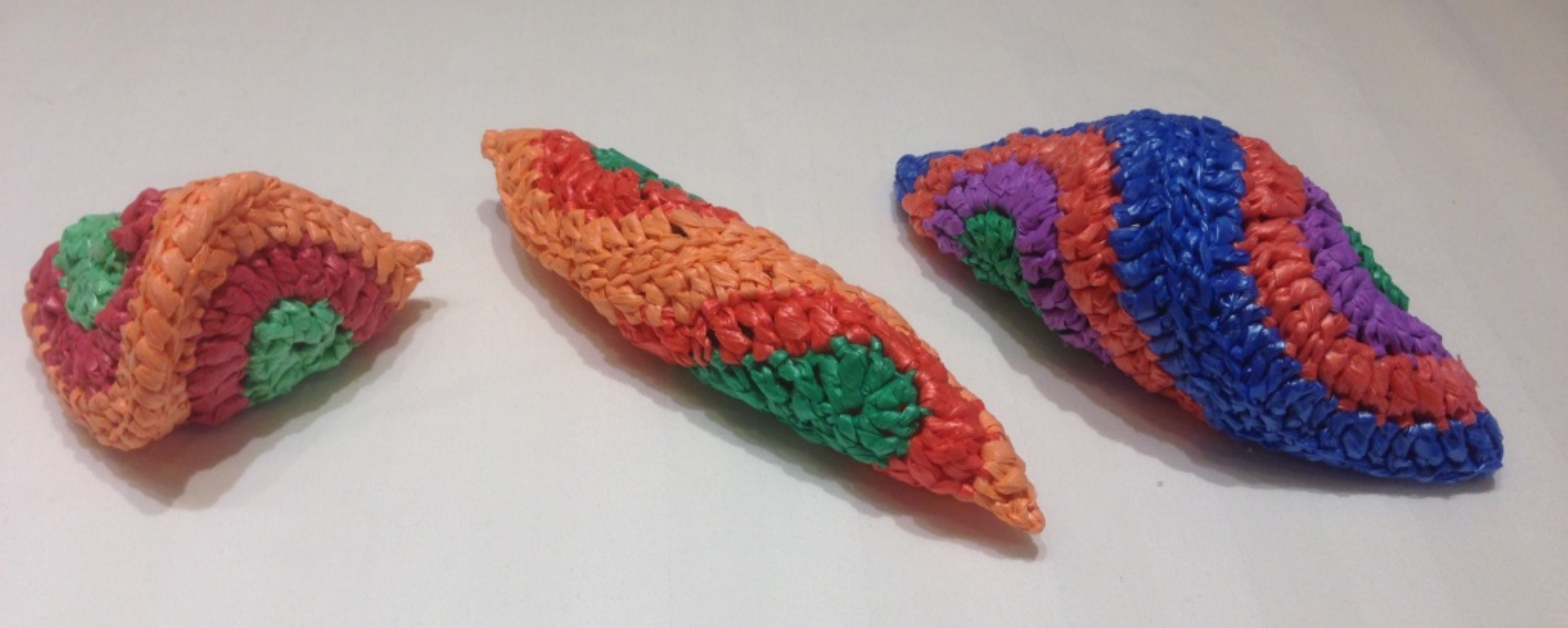}}}
\caption{\label{pita} Various pita-forms constructed from stadia by the author, by the method outlined in Appendix A. They each have the same value of $h$. The first two are identical apart from seam starting point. The third has a larger $r$ value. }
\end{center}
\end{figure}

But does The Shape belong on the sphericon family tree also?

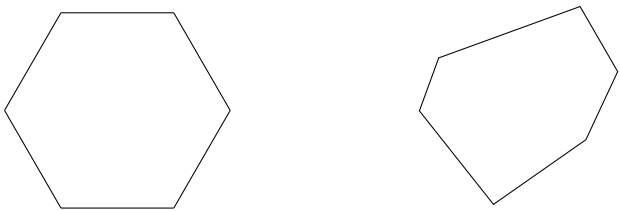
\begin{figure}
\begin{center}
\begin{tikzpicture}
  \begin{scope}
    \foreach \i in {0,...,5} {
      \coordinate (\i) at (\i*60:1.5);
    }
    \foreach \i in {0,...,5} {
      \draw let \n1={int(mod(\i+1,6))} in (\i) -- (\n1);
    }
  \end{scope}
  \begin{scope}[xshift=5cm,yshift=-1.25cm]
    \draw (0,0) -- ++(35:1.5) -- ++(65:1) -- ++(120:1) -- ++(200:2) -- ++(250:0.75) -- cycle;
  \end{scope}
\end{tikzpicture}
\caption{\label{hex} One is regular and one is irregular, but these are both hexagons.}
\end{center}
\end{figure}

As should have become clear from Section~\ref{solid}, any $N$-icon's surface is subject to symmetry, side length and angle constraints prescribed by its construction from a regular 3D object. What if one steps back from these origins, and hence from these constraints? Could not the stadium (Figure~\ref{stad}) be regarded as the outcome of relaxing  conditions which led to Figure \ref{temp}?  Just as the two shapes in Figure \ref{hex} are hexagons, one regular and one irregular, The Shape, with its one surface and one \textsf{S}-shaped seam between two distinct endpoints, could be termed an `irregular trisphericon', to coin a phrase.

To put this around the other way, the $(3,1)$-icon is a pita-form,  formed from a particular 2D shape (Figure \ref{temp}) which is {\em not} convex.

Moving on further from the original question about The Shape that motivated this investigation, what of other relationships between the $(N,k)$-icons, and pita- and D-forms? A kinship is foreshadowed in the works of the Turkish artist  Ilhan Koman. Consider the piece called {\em Rolling Lady} which dates from 1983 \cite{RL, AKA}. It consists of two bicones (four developable surfaces) interlocked at right angles in such a way as to roll. At its elegant heart lie two perpendicular circles, just as two perpendicular edges, circle sectors, are a striking feature of the sphericon. (In this article, the connection of the sphericon to the rollers and wobblers of \cite{EU}, and to the oloid, has not been explored, being somewhat tangential to the identification question under consideration.)

\begin{figure}
\begin{center}
{ \resizebox*{6cm}{!}{\includegraphics{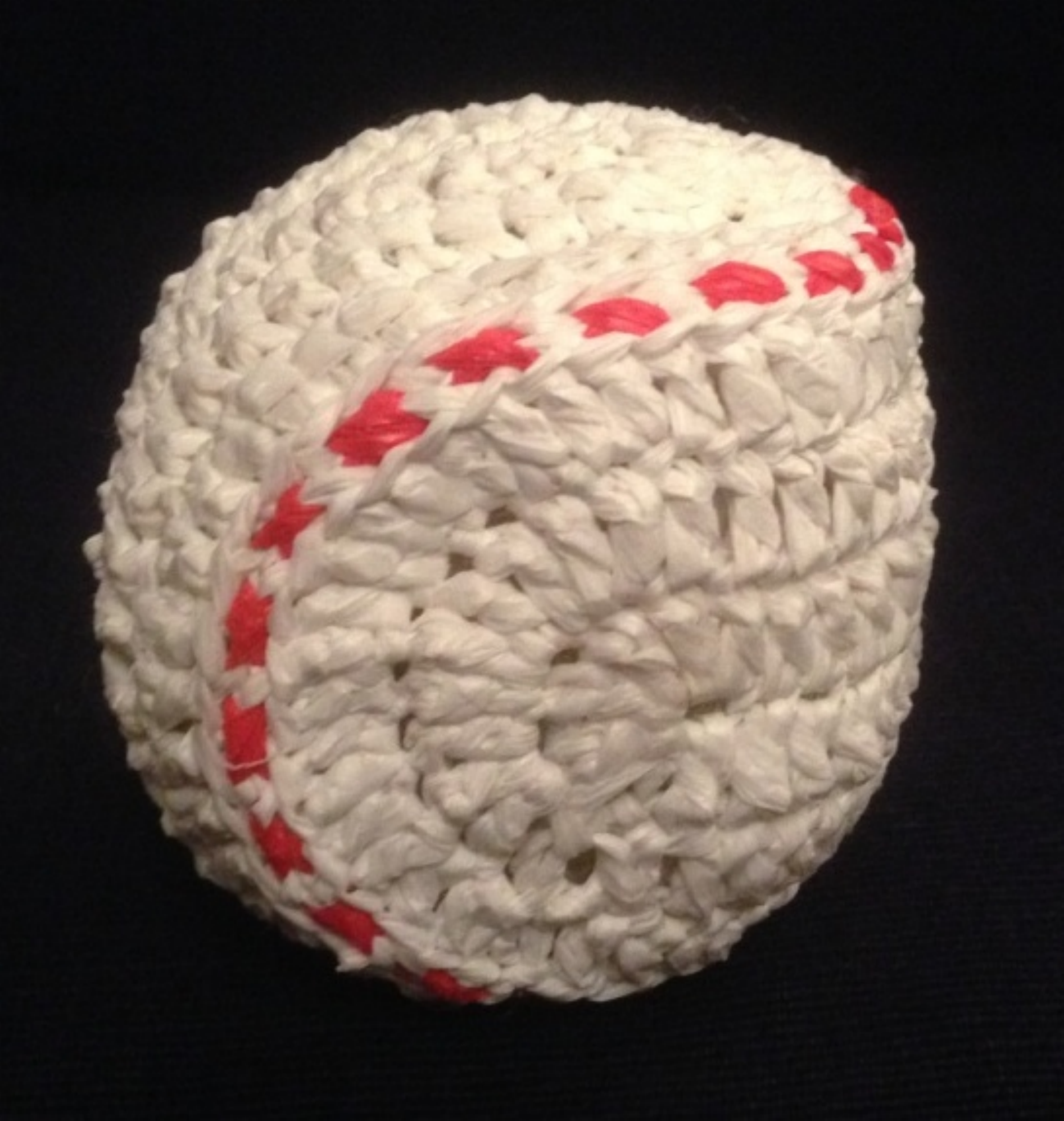}}}
\caption{\label{ball} {\em D-based}, a crochet figure by the author, February 2016. The title refers to the resemblance to a baseball highlighted by the colour choices, the insult (debasement) offered to this form in describing it as `uninteresting', and to the fact that it is, of course, a D-form.}
\end{center}
\end{figure}

Let's look at the shape that did not much interest Springett, the mm $(4,1)$-icon (Figure \ref{dual}). Again, the astute reader has probably pre-empted the statement that this is precisely a D-form obtained from two stadia, with $h= \upi r$, and attached to one another symmetrically. This inspired the crocheted piece {\em D-based} shown in Figure \ref{ball}.

D-forms made using stadia seem not to have been mentioned or visualised prior to 2007 \cite{DO}. Authors moved directly from dismissing the circle to trying out the ellipse. Perhaps this springs from the fact that the cone is probably the most well-known developable surface \cite{Sh1}, so that conic sections seemed a natural constituent for D-forms  (despite the fact that there is no formula for the perimeter of an ellipse, and matching perimeter length is the name of the game in making D-forms \cite{Bridgesnew}). Of course, a stadium has a perimeter that can be calculated as precisely as one might wish: $2h+2 \upi r$.

\begin{figure}
\begin{center}
{ \resizebox*{7cm}{!}{\includegraphics{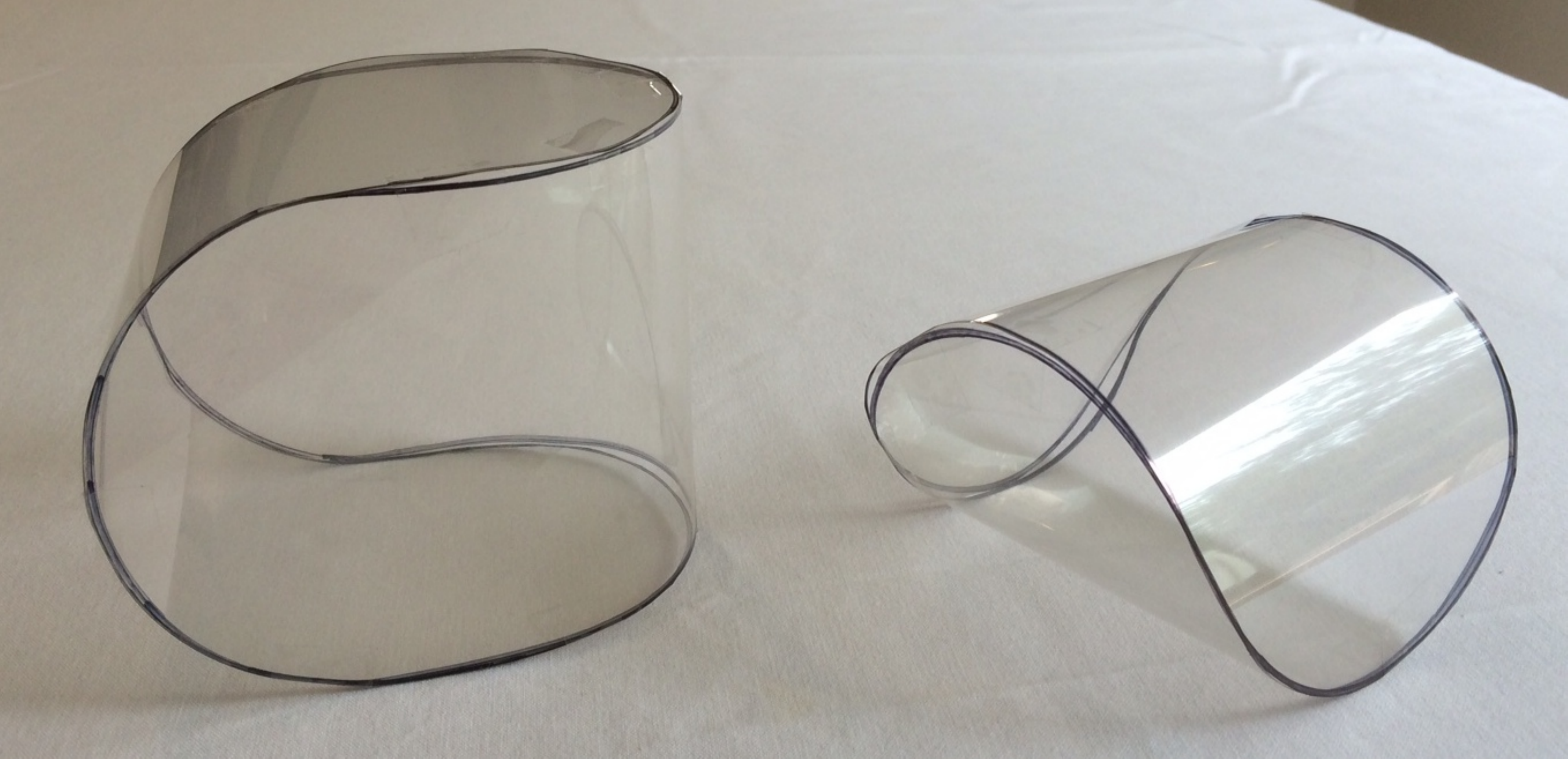}}}
\caption{\label{stable} {D-forms made from two stadia, attached to each other symmetrically. On the left, $h=4r$. On the right, $h=2r$.}}
\end{center}
\end{figure}

Two identical stadia can always be joined to form a D-form, but this will only be the mm $(4,1)$-icon if the midpoints of the semicircles on one stadium are joined to the midpoints of the straight sides of the other, and if $h=\upi r$. When $h>\upi r$, and the attachment is midpoint-to-midpoint, the elongated object is stable on one of its four semi-circular faces. Such a D-form makes an excellent four-sided die. In fact, it is arguably nicer than the one based on the tetrahedron, which has a point, not a face, in the air! When $h< \upi r$, and the attachment is symmetric, the object is stable only when it rests on the centre-line of one or other of its two curved surfaces. If tilted slightly from this position, it will rock back.  Hollow mock-ups of these two D-forms can be seen in Figure \ref{stable}. A perfectly constructed mm $(4,1)$-icon has stability properties that transition between these two cases.
Tilting the mm $(4,1)$-icon when it is sitting on the centre-line of one of its two curved surfaces does not raise or lower the centre of gravity, so that this is a position of neutral equilibrium; it rolls, until it sits flat on one of its two semicircular faces.

\section{Generalising}\label{concl}

It is interesting to classify which other $(N,k)$-icons are  D-forms and pita-forms (in the original sense, constructed using 2D shapes with no particular convexity requirement imposed). The information needed to perform this classification is available, though using unconventional mathematical nomenclature, and at least one error, in early studies at un-maintained web sites \cite{PR, RK}. It is collected and corrected in Table \ref{classify}. In interpreting this table, the values of $k$ are restricted to $0\leqslant k <\frac{N}{2}$ for $N$ odd, and to $0\leqslant k \leqslant\frac{N}{4}$ for $N$ even, limiting the forms to those which are unique (up to chirality).

\begin{table}
\tbl{Features of $(N,k)$-icons}
{\begin{tabular}[l]{cccc}
\toprule
   & $N$ odd & vv $N$ even
         & mm $N$ even
          \\

\colrule
  Closed edges &$\frac{\text{gcd}(N, k)-1}{2}$ & $\text{gcd}(\frac{N}{2},k)-1$&$\text{gcd}(\frac{N}{2},k)$  \\
  &&&\\
  Edges running between two distinct endpoints & 1 & 2 & 0  \\
  &&&\\
 Surfaces that can be `traced' & $\frac{\text{gcd}(N, k)-1}{2}$&$\text{gcd}(\frac{N}{2},k)$ & $\text{gcd}(\frac{N}{2},k)-1$\\
 &&&\\
 Other surfaces &1&0&2\\
\botrule
\end{tabular}}
\label{classify}
\end{table}

An $(N,k)$-icon will be a pita-form if it has one surface and one edge that runs between two distinct endpoints; it will be a $D$-form if it has two surfaces and one closed edge.

From this we can conclude:
\begin{itemize}
\item When $N$ is odd, and co-prime to $k$, the $(N,k)$-icon is a generalised pita-form. (This includes the trisphericon. The first values for which this does not apply are $N=9,\  k=3$.)
\item When $N$ is even, and the $(N,k)$-icon has been generated using the midpoint-to-midpoint symmetry axis, it will be a D-form when $\frac{N}{2}$ and $k$ are co-prime. (This includes, for example, all 5 possible $N=22$ forms, but only 2 of the possible 6 $N=12$ forms. )
\end{itemize}

\begin{figure}
\begin{center}
{ \resizebox*{6cm}{!}{\includegraphics{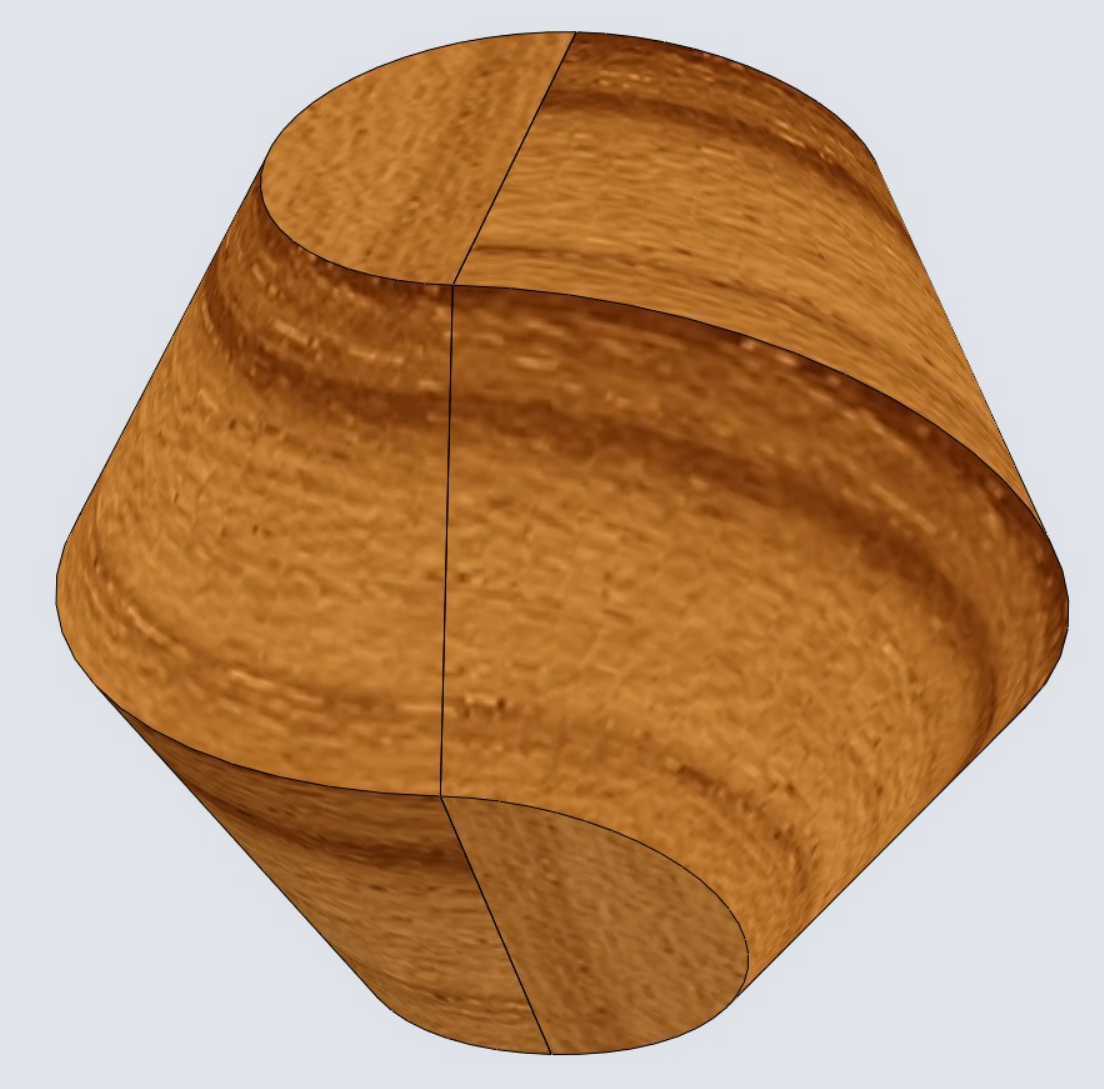}}}
\caption{\label{hexasphericon} A wooden hexasphericon generated using the midpoint-to-midpoint symmetry axis. Permission given to reproduce from \cite{SS}.}
\end{center}
\end{figure}

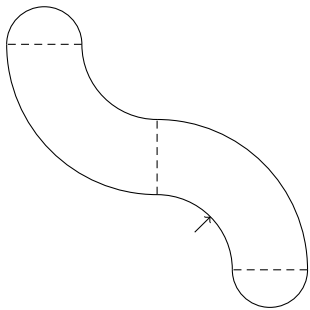
\begin{figure}
 \begin{center}
\begin{tikzpicture}
  \draw (0,0) arc [start angle=0, end angle =90, radius=2];
\draw (0,0) arc [start angle=0, end angle=-180, radius=0.5];
\draw (-1,0) arc [start angle =0, end angle =90, radius =1];
\draw (-2,1) arc [start angle = -90, end angle =-180, radius =2];
\draw (-2,2) arc [start angle =-90, end angle =-180, radius = 1];
\draw (-3,3) arc [start angle = 0, end angle =180, radius =0.5];
\draw[->] (-1.5,0.5) -- (-1.29,0.707);
  \draw[densely dashed] (0,0) --  (-1,0);
  \draw[densely dashed] (-3,3) --  (-4,3);
    \draw[densely dashed] (-2,1) --  (-2,2);

  \end{tikzpicture}
  \end{center}
\caption{\label{hextemp} The template for each of the two surfaces that comprise the mm $(6,1)$-icon. The smaller arcs have the same length as the semi-circles, and the larger arcs are twice as long. The arrow indicates where the midpoint of one semicircle of the other surface should be attached to begin the D-forming. }
\end{figure}

By way of illustration, Figure~\ref{hexasphericon} shows a wooden version of the mm $(6,1)$-icon (one possible hexasphericon); the two semicircular ends that are visible are from the same surface. When a regular hexagon (side length $r$) is rotated to form a solid of revolution (about the midpoints of a pair of opposite edges), its vertices trace out three circles. Two of these have diameter $r$, and the other has diameter $2r$, the distance between opposite vertices in a regular hexagon, by elementary geometry (see Figure~\ref{hex}). Each of these circles is cut in half, relocated in space and then joined to half of another circle, when the cut-rotate-reglue process is carried out. Apart from the circular ends, the other parts of this  hexasphericon's surface are pieces of cones, and each is formed (developed) from a quarter of an annulus. Thus each of the two surfaces that form the mm $(6,1)$-icon by D-forming is of the shape shown in Figure~\ref{hextemp}.

\begin{figure}
 \begin{center}
\begin{tikzpicture}

 \draw (0,0) --  (3.8,0);
  \draw (0,1) --  (3.8,1);
  \draw (3.8,1) arc [start angle=90, end angle =-37.3, radius = 1.71];
  \draw (3.8,0) arc [start angle=90, end angle =-37.3, radius = 0.71];
  \draw (0,0) arc [start angle =-90, end angle = -217.3, radius =1.71];
  \draw (0,1) arc [start angle =-90, end angle= -217.3, radius =0.71];
  \draw (-0.56, 2.14) arc [start angle = -37.3, end angle = 142.7, radius =0.5];
  \draw (4.36, -1.14) arc [start angle = -217.3, end angle = -37.3, radius =0.5];

  \end{tikzpicture}
  \end{center}
\caption{\label{dualoct} The template for each of the two surfaces that comprise the mm $(8,1)$-icon, one of the possible octasphericons. }
\end{figure}
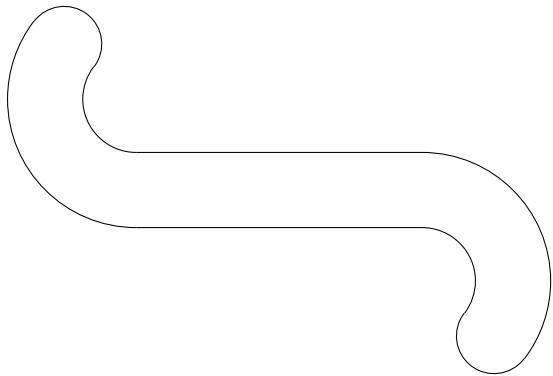

For the mm $(8,1)$-icon (an octasphericon) a similar geometric argument shows that each of its constituent surfaces is made up of a rectangle, two semi-circles and two sectors (opening angle ${\upi}/{\sqrt{2}}$) from an annulus. The inner arc of the annulus sector has the same length as the semi-circles ($\upi r$), and the length of the outer arc and of the rectangle side is $(1+\sqrt{2})\upi r$. The template is shown in Figure~\ref{dualoct}.

For $N=10$, there are two possible mm $(10,k)$-icons. Their two constituent surfaces are each bounded by 2 semicircles, 2 arcs of length $\upi r$ and $(1+2\sin(\frac{\pi}{5})+2\sin(\frac{\pi}{10}))\upi r$, and four of length $(1+2\sin(\frac{\pi}{5}))\upi r$. The order in which they connect to each other is different for $k=1$ and $k=2$.

In the static representation of words and diagrams, it is not easy to appreciate the features of these forms. As well as making some of their own, readers are encouraged to enjoy some remarkable on-line animations of them (for numerous higher values of $N$ and $k$) \cite{RK}.

As far as can be judged from the published literature, these instructions for creating the `piece-wise circular' surfaces of the mm $(N,k)$-icons by D-forming have not been explicitly stated previously, though presentations on sphericons and D-forms have been made at the same conference \cite{AKA,Sp,Wi}! This may be because of the very different geography of their birthplaces, in the realm of solids and of surfaces respectively. It seems to have taken crochet fibre to link them.

\section*{Acknowledgements} I wish to acknowledge the contributors to the thread {\em How would you call this shape?} in the Woolly Thoughts forum of Ravelry, in particular Annagret,  bevbh, irishlacenet, MadeleineS, soyloquesoy and Woolhelmina, for their inspiration of this piece, and for pooling their collective wisdom about fibre arts and mathematics. Thanks also to Jane Pitkethly for creating four of the line drawings and to the fibre artists who permitted photos of their work to be used. Particular thanks to the estate of the late Steve Mathias for giving permission for photos of his wood-turning to be used, and to Craig Kaplan for Figure 3. The perceptive comments of editor and associate editor have also improved this paper significantly, and I thank them.

\medskip

\appendices
\section{Make your own dual sphericon or pita-form}
To crochet a stadium is not as daunting a task as making the Lorentz Manifold \cite{HO} or even a hyperbolic plane \cite{DT}. To make a dual sphericon, you will need to crochet two identical pieces. To make a pita-form, you need only one. You will need yarn, a crochet hook, a sewing needle and scissors (only to cut threads, not to cut the surface!).

Any (free or purchased) non-lacy pattern for an `oval' will suffice. This could be a pattern for a rug or table mat. `Oval' is used in craft forums (and daily conversation) for just about anything  which is not symmetric enough to be a circle  and not  angular enough to be a rectangle, and even for 3D forms, so check to see that it is actually going to produce a stadium (Figure \ref{stad}). If you are familiar with both knitting and crochet, you will realise that to produce a stadium, crochet is easiest. This is related to the reason that crochet is more natural for producing hyperbolic surfaces \cite{DT}; because one works with a stitch at a time on the hook, and can control stitch height, different types of shapes can be made in crochet from those made with knitting.

I have chosen to make my models so far using plarn (plastic yarn, made by recycling shopping bags) which produces quite a stiff fabric even for treble (UK) stitches. Using acrylic yarn and a smaller-than-usual hook to produce a dense fabric in double (UK) crochet would also work \cite{DT}. If you decide to stuff your models, use just enough stuffing that they do not collapse, but not so much as to distort them. My models are quite small, but larger models will probably need stuffing. Another method to get around the inherent stretchiness and floppiness of most knitted or crocheted fabric, and the shaping disadvantages of knitting, is to use felted knitted fabric, cut to shape, as seen in Figure~\ref{MS}. I closed my models by overstitching, pulling only tight enough to close the seam, again to avoid distortion.

For the dual sphericon, begin by joining the midpoint of a semicircle on one stadium to the midpoint of the straight side of the other, and then join stitch by stitch, until the two pieces are joined. For the pita-form, you can pick any starting point (apart from the exact midpoints) for your seam. It's fun to make several, using different starting points.

To control the size of the rectangular part of your stadium,  vary the length of the starting chain. The diameter of the semi-circular parts depends on the stitch height (double or treble) and the number of rounds completed; you will need to make a gauge swatch if you want stadia of precise measurements. It might be fun to use $h/(\upi r)=$

\begin{align*}
1+2(1+\sqrt{2})=3+2\sqrt{2}&\approx 5.8284\\
4+6\sin(\tfrac{\pi}{5})+2\sin(\tfrac{\pi}{10})=4+3 \sqrt{(5-\sqrt{5})/2}+(\sqrt{5}-1)/2&\approx 8.1447.
\end{align*}
which correspond to relaxing the templates of Figure~\ref{hextemp} and Figure~\ref{dualoct} into stadia; the objects produced could be called irregular dual $(N,k)$-icons!

\begin{figure}
\begin{center}
{\resizebox*{6cm}{!}{\includegraphics{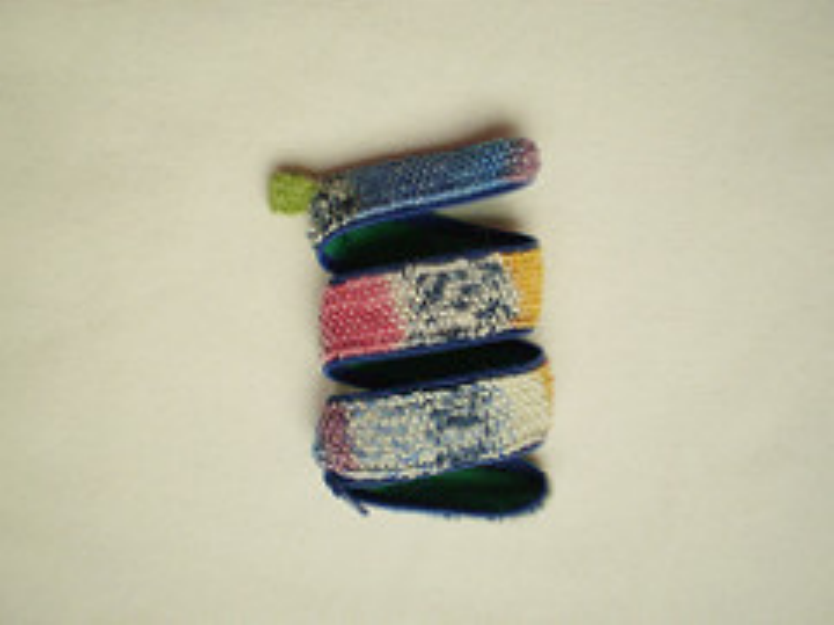}}}
\caption{\label{zippy} {\em Zippy Strip}. Courtesy of Frankie Brown.}
\end{center}
\end{figure}

\begin{figure}
\begin{center}
{\resizebox*{3.5cm}{!}{\includegraphics{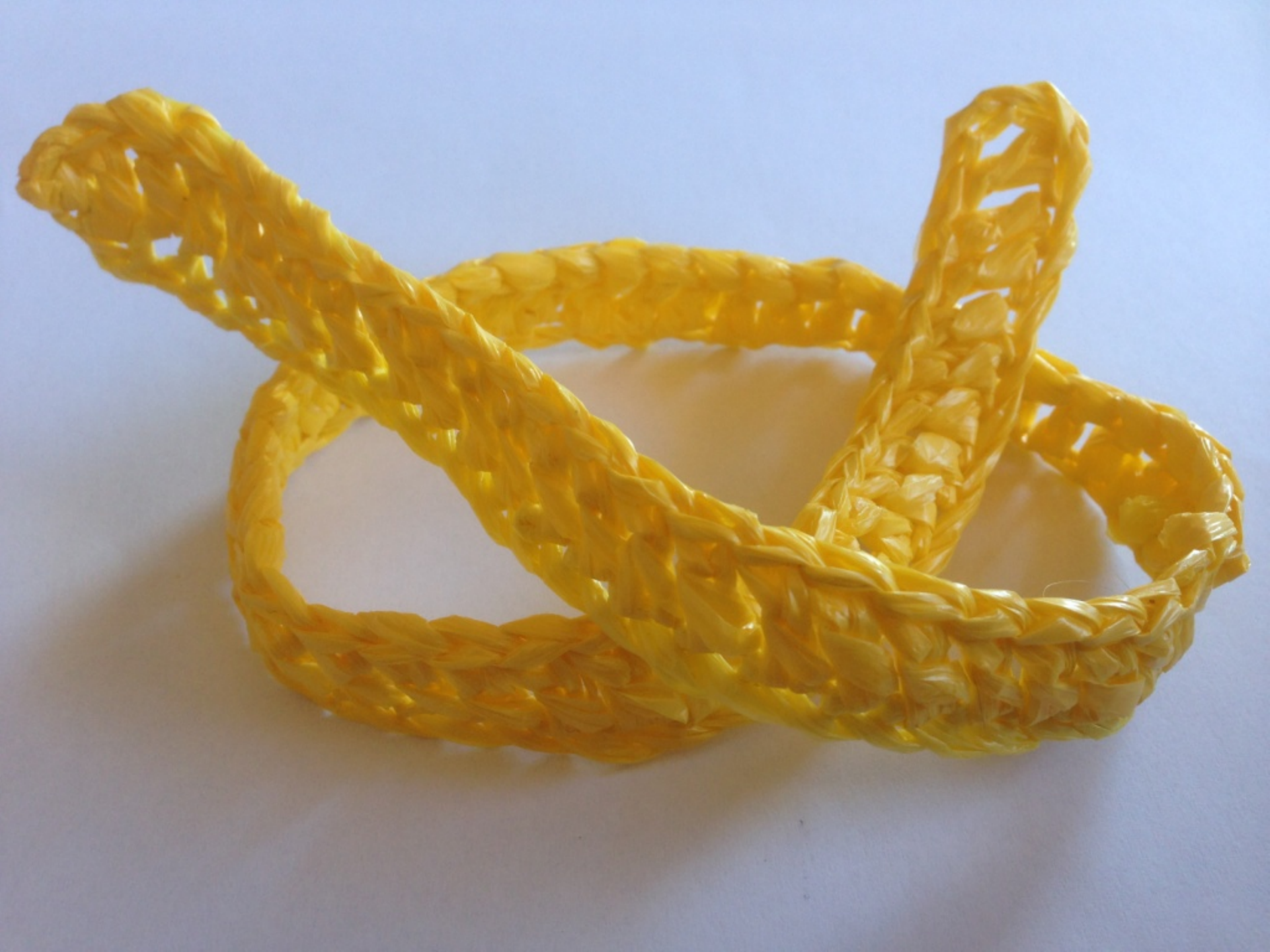}}}
{\resizebox*{3.5cm}{!}{\includegraphics{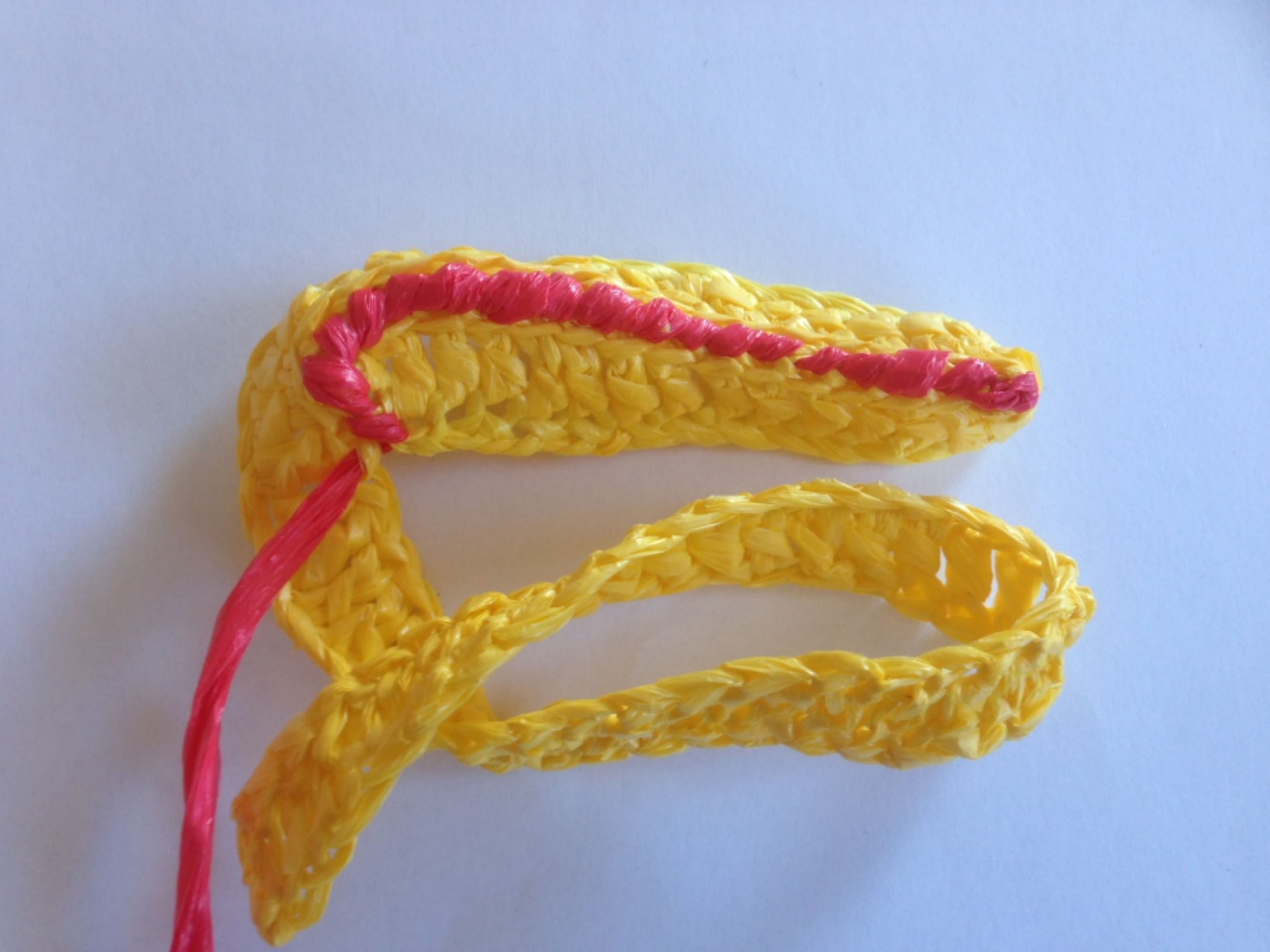}}}
{\resizebox*{3.5cm}{!}{\includegraphics{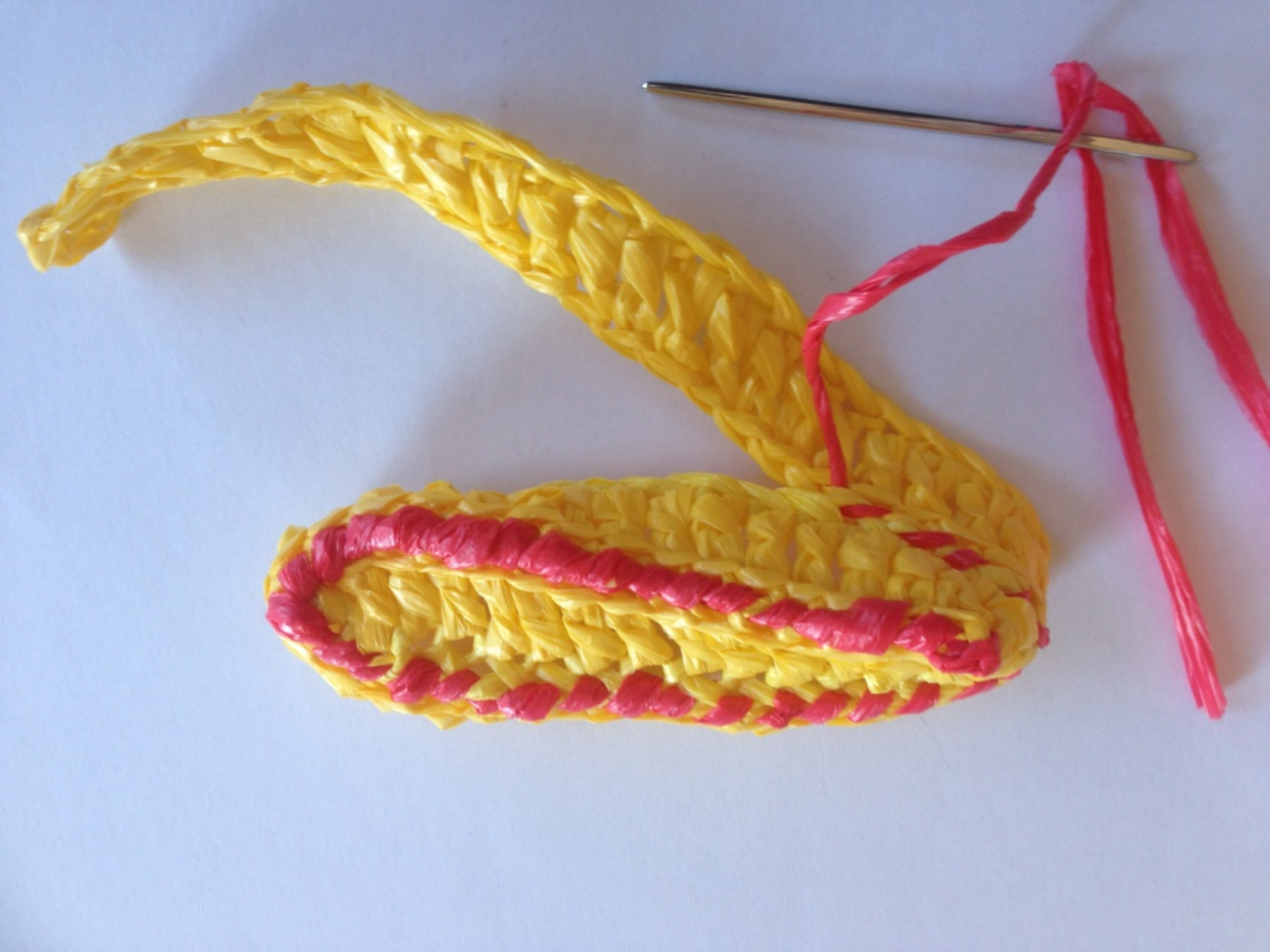}}}
{\resizebox*{3.5cm}{!}{\includegraphics{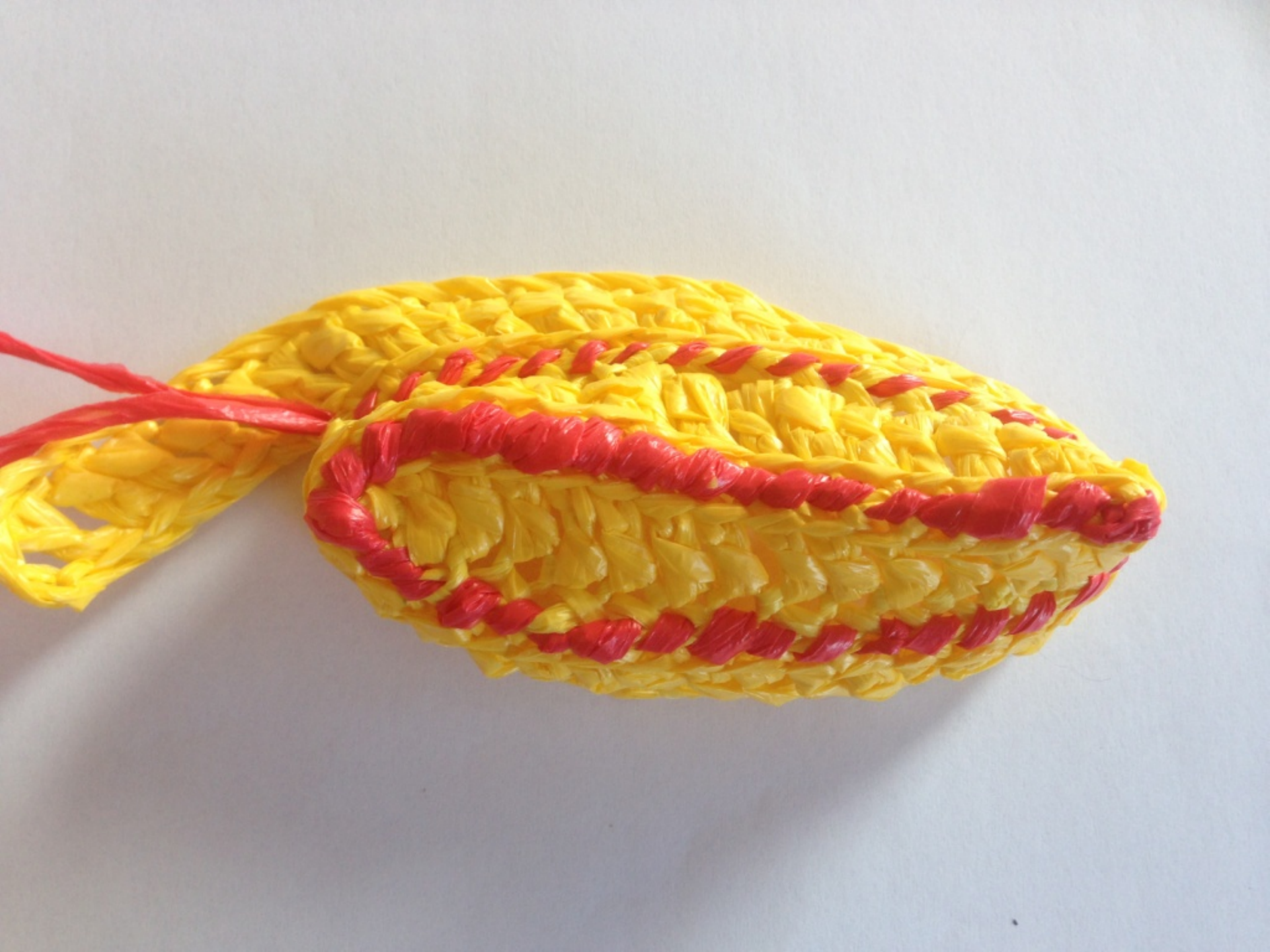}}}
{\resizebox*{3.5cm}{!}{\includegraphics{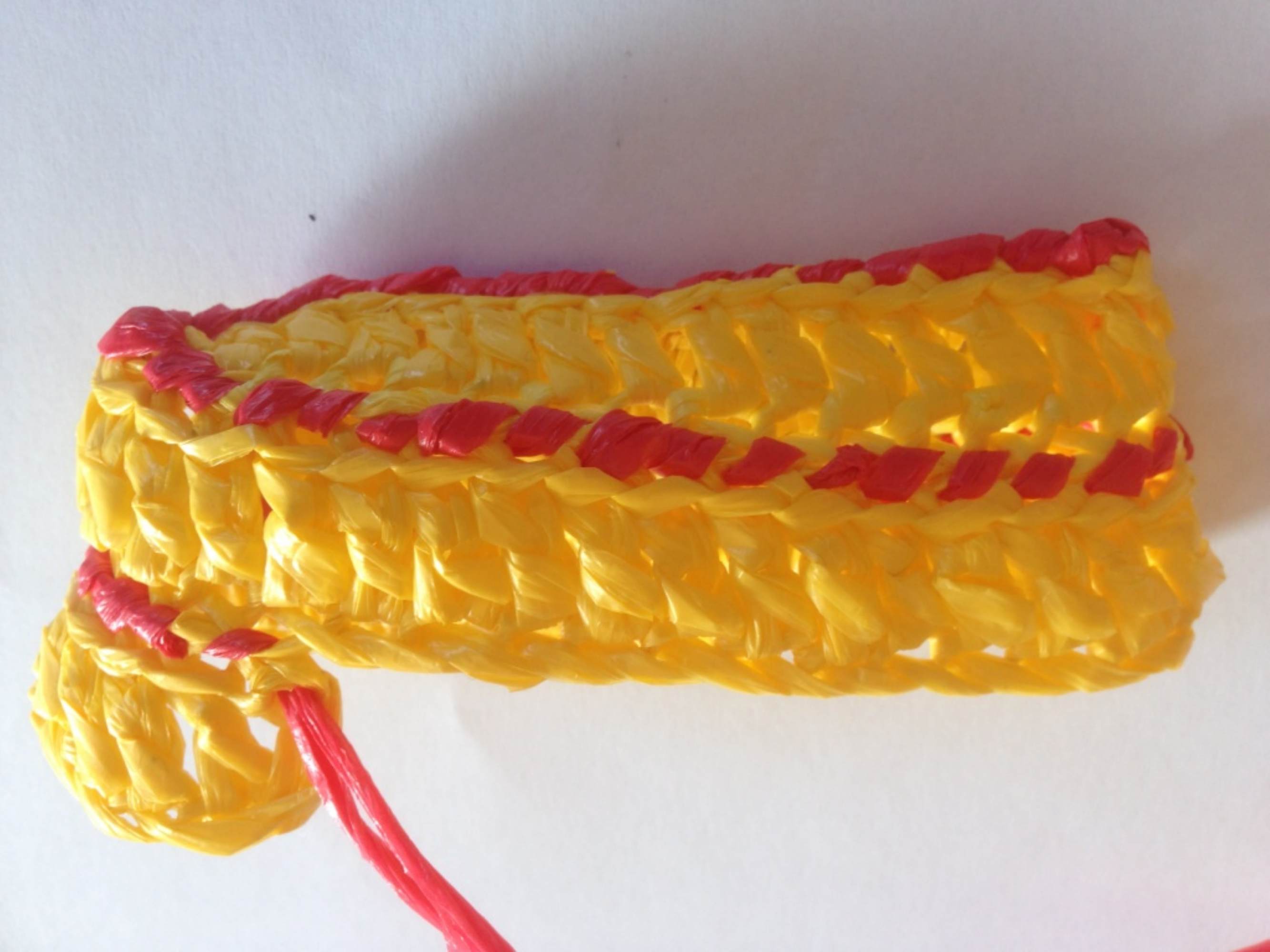}}}
{\resizebox*{3.5cm}{!}{\includegraphics{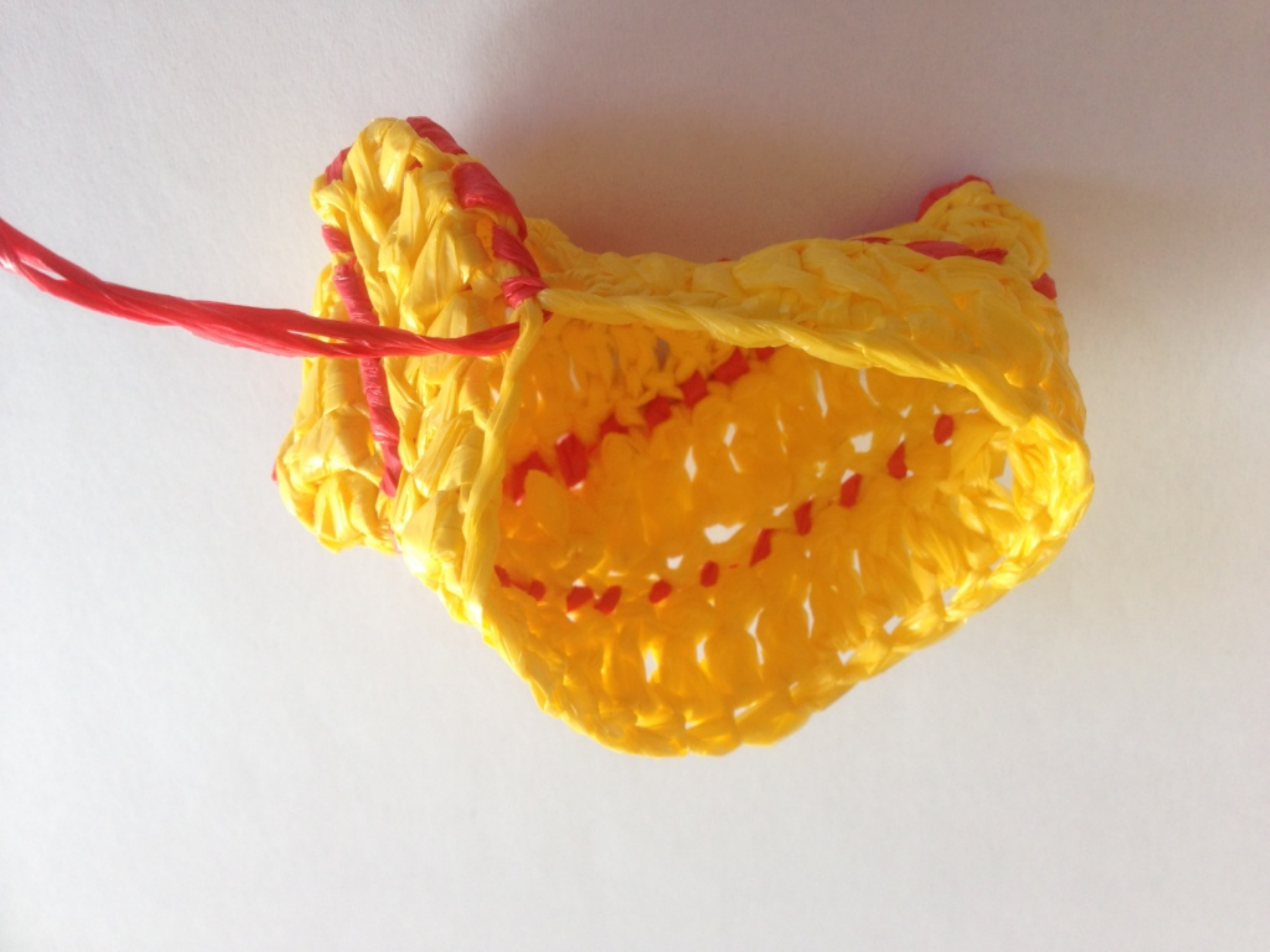}}}
{\resizebox*{3.5cm}{!}{\includegraphics{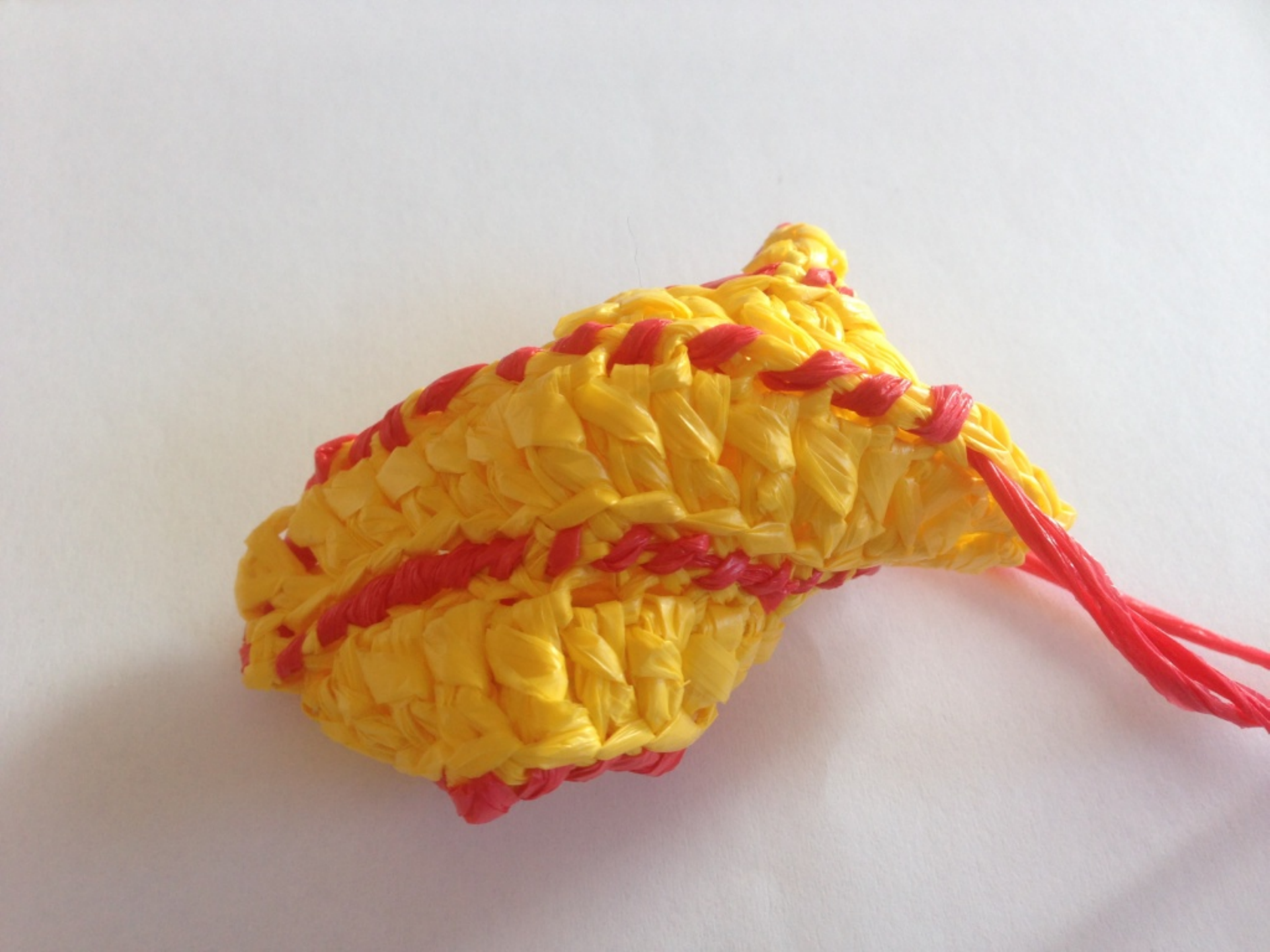}}}
{\resizebox*{3.5cm}{!}{\includegraphics{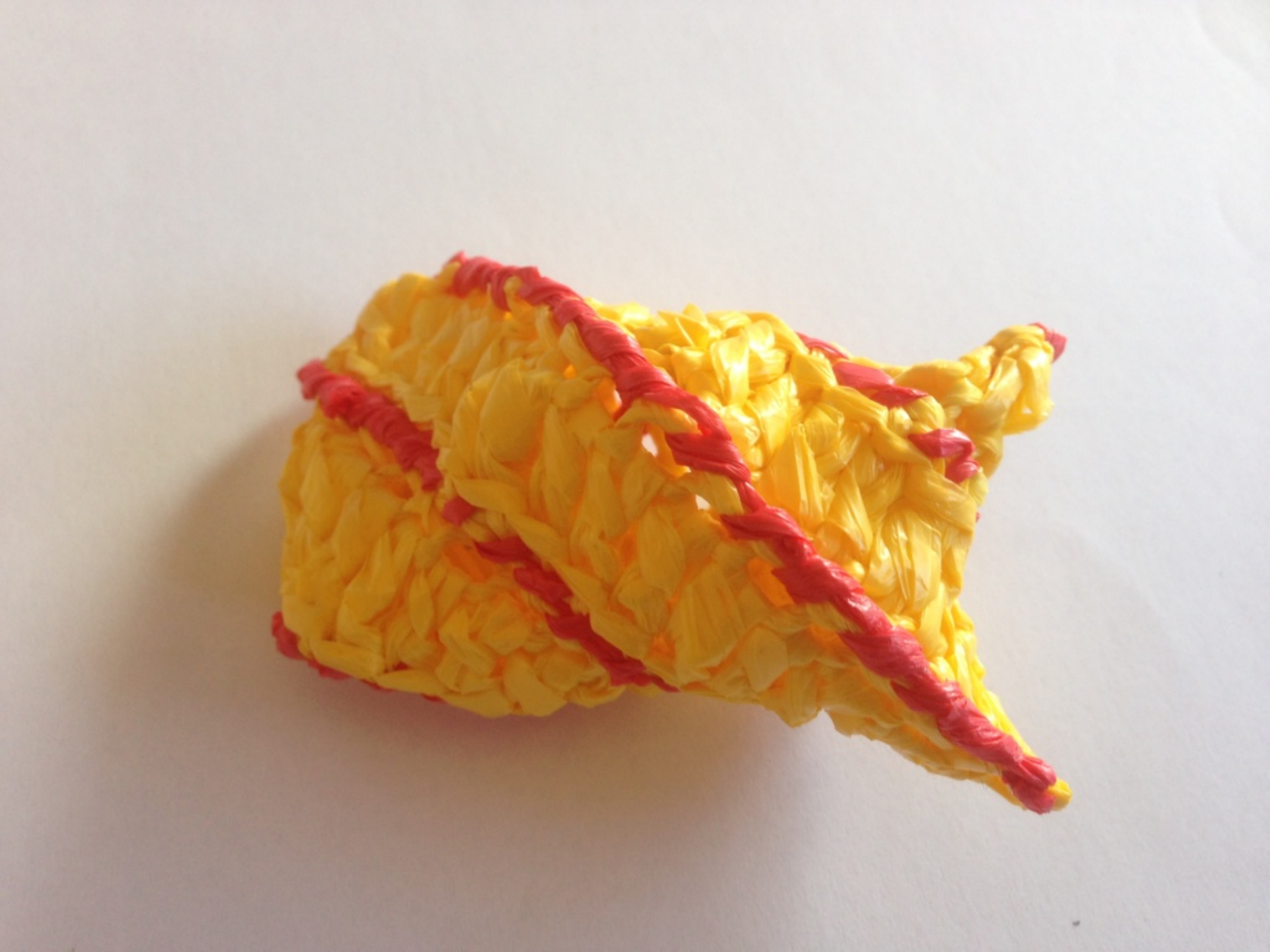}}}
\caption{\label{contort} A stadium of radius $r=1.5$ cm and $h=50$ cm is seamed to itself. While it begins as does the Zippy Strip, it develops a different shape as the final part of the seam is made.}
\end{center}
\end{figure}

Interesting pita-forms are obtained for $h \gg r$ such as those created using the pattern {\em Zippy Strip} \cite{ZS}. In this case, the stadium continues to wrap around itself, with the seam forming a helix in space (rather than just an \textsf{S}). Again, the starting point of the seam will determine the final outcome. Some makers describe a flattish purse (Figure \ref{zippy}), while others obtain something more `pyramidal'. In making one myself, I observed that remarkable feature of D-forms, that they continue to develop right up until they are closed. This development is captured in Figure \ref{contort}.

\label{lastpage}

\end{document}